\documentclass[11pt]{amsart}
\usepackage{amsmath,amscd,amsfonts,amssymb,latexsym,amsthm,array}

\usepackage{xcolor}
\usepackage{enumitem}
\theoremstyle{plain}
\newtheorem{thm}{Theorem}
\newtheorem{prop}[thm]{Proposition}
\newtheorem{lemma}[thm]{Lemma}

\newtheorem{cor}[thm]{Corollary}

\newtheorem{remark}[thm]{Remark}

\theoremstyle{definition}

\theoremstyle{remark}
\newtheorem{rem}[thm]{Remark}

\DeclareMathOperator{\spn}{span}
\DeclareMathOperator{\supp}{supp}
\DeclareMathOperator{\sgn}{sgn}

\newcommand{\N}{\mathbb{N}}
\makeatletter
\@namedef{subjclassname@2020}{%
  \textup{2020} Mathematics Subject Classification}
\makeatother

\allowdisplaybreaks

\begin{document}

\title{On uniqueness and plentitude of subsymmetric sequences}
%\dedicatory{}

%------------------------------------------------------------------------------------------------------------------------------------------------

\author[P. G. Casazza]{Peter G. Casazza}

\address{Department of Mathematics, University of Missouri, Columbia, MO 65211-4100, USA}

\email{casazzap@missouri.edu}

%------------------------------------------------------------------------------------------------------------------------------------------------

\author[S. J. Dilworth]{Stephen J. Dilworth}

\address{Department of Mathematics, University of South Carolina, Columbia,
SC 29208, USA}

\email{dilworth@math.sc.edu}

%------------------------------------------------------------------------------------------------------------------------------------------------

\author[D. Kutzarova]{Denka Kutzarova}

\address{Department of Mathematics, University of Illinois Urbana-Champaign,
Urbana, IL 61807, USA; Institute of Mathematics and Informatics, Bulgarian Academy of Sciences, Sofia, Bulgaria}

\email{denka@illinois.edu}

%------------------------------------------------------------------------------------------------------------------------------------------------

\author[P. Motakis]{Pavlos Motakis}

\address{Department of Mathematics and Statistics, York University, 4700 Keele Street, Toronto, Ontario, M3J 1P3, Canada}

\email{pmotakis@yorku.ca}

%------------------------------------------------------------------------------------------------------------------------------------------------

\thanks{The first author was supported by
 NSF DMS 1906025.  The second author was supported by Simons Foundation Collaboration Grant No. 849142.
 The third author was supported by Simons Foundation Collaboration Grant No. 636954. The fourth author was supported by NSERC Grant RGPIN-2021-03639.}

%\date{\today}

%\keywords{}


\begin{abstract}
We explore the diversity of subsymmetric basic sequences in spaces with a subsymmetric basis. We prove that the subsymmetrization $Su(T^*)$ of Tsirelson's original Banach space provides the first known example of a space with a unique subsymmetric basic sequence that is additionally non-symmetric. Contrastingly, we provide a criterion for a space with a subsymmetric basis to contain a continuum of nonequivalent subsymmetric basic sequences and apply it to  $Su(T^*)^*$. Finally, we provide a criterion for a subsymmetric sequence to be equivalent to the unit vector basis of some $\ell_p$ or $c_0$.
\end{abstract}

\maketitle

\section{Introduction}
A main question of the structure theory of Banach
spaces is whether any infinite-dimensional space contains an infinite-dimensional subspace which
is isomorphic to a space from a list of spaces with ``nice''
properties. The most natural first
question was if any Banach space contained an isomorphic copy of $c_0$
or $\ell_p$, $ 1\le p < \infty$, or more generally, a symmetric basic
sequence. Recall that a sequence $(x_j)_{j=1}^{\infty}$ is a basic
sequence if it is a (Schauder) basis of its closed linear span; two
basic sequences $(x_j)_{j=1}^{\infty}$ and $(y_j)_{j=1}^{\infty}$ are
said to be equivalent provided a series
$\sum_{j=1}^{\infty}a_{j}x_{j}$ converges if and only if
$\sum_{j=1}^{\infty}a_{j}y_{j}$ does.  A basic sequence
$(x_j)_{j=1}^\infty$ is \emph{symmetric} if the rearranged sequence
$(x_{\pi(j)})_{j=1}^{\infty}$ is equivalent to $(x_j)_{j=1}^\infty$
for any permutation $\pi$ of $\mathbb{N}$. The above question was solved in the
negative by Tsirelson in 1974 \cite{T} and his space led to remarkable new
developments. Actually, what is now
referred to as the Tsirelson space $T$ , is the construction given by Figiel and Johnson \cite{FJ} and
the original Tsirelson space is its dual $T^*$ .

The class of \emph{subsymmetric} basic sequences, i.e., unconditional and equivalent to all of
their subsequences, is formally more general than the class of
symmetric ones. For a while, these two concepts were believed to be
equivalent until Garling \cite{G} provided a counterexample.
However, subsymmetric bases, far from being just a natural 
generalization of symmetric bases, later played an important role by themselves within
the general theory. The first non-arbitrarily distortable space
constructed by Schlumprecht \cite{S} has a subsymmetric basis.

Altshuler \cite{A} (see also  \cite{LT1}) constructed a space in which all symmetric basic sequences
are equivalent to its symmetric basis. Later, another example of this
kind was
built in \cite{CS}, based on a symmetric version $S(T^*)$ of $T^*$. Recently, Albiac, Ansorena and
Wallis  \cite{AAW} used Garling-type spaces
 to provide the first example of a Banach space with a unique subsymmetric
basis which is not symmetric. However, as shown in a sequel paper \cite{AADK}, that space
contains a
continuum of non-equivalent subsymmetric basic sequences. A careful inspection of Altshuler's proof shows that, in fact,
all subsymmetric basic sequences of his space are also equivalent to
the symmetric basis. The same turned out to be true in the case of
$S(T^*)$.

In view of Altshuler's example, it was asked in
\cite{KMP} and \cite{AADK} whether there exists
a space with a subsymmetric basis with a unique, up to equivalence,
subsymmetric basic sequence which is not symmetric.
In this paper we
answer this question in the positive. For this purpose, we use the
subsymmetric version $Su(T^*)$ \cite{CS} of the
original Tsirelson space $T^*$ \cite{T}. One can
also define a subsymmetric version of Altshuler's space, however we do
not know whether all its subsymmetric basic sequences are
equivalent. In addition, we show that the dual space $Su(T^*)^*$ has
no symmetric basic sequences. Next,  we extend to subsymmetric bases Altshuler's characterization \cite{A1} of the unit vector bases of $c_0$ and $\ell_p$
among symmetric bases and
we give various sufficient conditions
for a space with a subsymmetric (resp. symmetric) basis to have a continuum of
subsymmetric (resp. symmetric) basic sequences. These results are applied to $Su(T^*)^*$ and to the dual of Altshuler's space.
In the final section we give some applications to the  set of  spreading models of a given Banach space.

\section{A space with a unique subsymmetric basic sequence}
We begin with some definitions. Given two basic sequences
$(x_n)_{n=1}^\infty$ and $(y_n)_{n=1}^\infty$ in Banach spaces $X$ and
$Y$, respectively, we say that $(x_n)_{n=1}^\infty$ \emph{$K$-dominates}
$(y_n)_{n=1}^\infty$ if there is a bounded linear operator
$T\colon [(x_n)_{n=1}^\infty] \to [(y_n)_{n=1}^\infty]$, $\|T\| \le K$, such that
  $T(x_n)=y_n$ for all natural numbers $n$. We say that  $(x_n)_{n=1}^\infty$ \emph{dominates}
  $(y_n)_{n=1}^\infty$, denoted $(y_n)\preceq (x_n)$,  if  $(x_n)_{n=1}^\infty$ $K$-dominates
$(y_n)_{n=1}^\infty$ for some $K$.
  A \emph{block basis} with
respect to a basic sequence $(x_n)_{n=1}^\infty$ is a sequence
$(y_n)_{n=1}^\infty$ of non-zero vectors of the form
$y_n=\sum_{k=p_n+1}^{p_{n+1}} a_k x_k$ where $p_1 < p_2 < \dots$ is an
increasing sequence of natural numbers.
For a vector $x$ in the closed linear span of $(x_n)_{n=1}^\infty$,
its support (with respect to $(x_n)_{n=1}^\infty$) is the set of
indices of its non-zero coefficients. For finite sets of natural
numbers $E$ and $F$ we say that $E < F$  if $\max(E)
<\min(F)$.
For a natural number $n$, we say $n < x$, resp. $n \le x$, if $n <
\min(\supp(x))$, resp. $n \le  \min(\supp(x))$. A basic sequence $(x_n)$ is called \emph{1-subsymmetric} if it is $1$-unconditional and isometrically equivalent to its subsequences.

Recall the Figiel-Johnson construction of $T$. For each $x=\sum_n
a_nt_n$, where $(t_n)$ is the canonical basis of $T$,
$$
\|x\| = \max\Big\{ \max_n |a_n|, \frac12 \sup \sum_{j=1}^k \|E_j
  x\|\Big\},
$$
where the inner supremum is taken over all choices $k \le E_1 < E_2<
\dots <E_k$ and for a set of natural numbers $E$, $Ex$ is the natural
basis projection onto $spn\{t_j \colon j\in E\}$.

From the definition of $T^*$ it follows that if $x_1 < x_2 < \dots
x_k$ is a normalized block basis of $(t_n^*)$ with $k \le x_1$, then
$\left\| \sum_{j=1}^k x_j\right\| \le 2$.

The subsymmetric version of $T^*$ \cite{CS} denoted $Su(T^*)$ has the norm:
$$
\left\|\sum a_nt_n^*\right\|=\sup\left\{\left\|\sum a_{n_i}t^*_i\right\|_{T^*} \colon
n_1<n_2<\dots \right\}.
$$
Although this is not entirely ovious, $\|\cdot\|$ is subsymmetric (Corollary \ref{subsymmetric norm}).

\begin{prop}\label{pr1}
The canonical basis $(t^*_n)$ of $Su(T^*)$ is not symmetric.
\end{prop}

This was remarked in \cite{CS}. We shall prove a stronger result later
in this paper (see Section \ref{dual empty of symmetry}).

\begin{thm}\label{th2}
Every subsymmetric sequence in $Su(T^*)$ is equivalent to $(t_i^*)$.
\end{thm}

Our proof will be given in several steps and it uses the following facts about $T^*$ from \cite {CS}.

\begin{prop}
\label{pr2}~
\begin{enumerate}[leftmargin=19pt,label=(\arabic*)]

\item\label{equivalence by range} There exists a constant $K$ such that if $x_n=\sum_{p_n+1}^{p_{n+1}} a_it^*_i$
and $y_n=\sum_{p_n+1}^{p_{n+1}} b_it^*_i$ ($n\ge 1$) are normalized
block bases then $(x_n)$ and $(y_n)$ are $K$-equivalent.

\item\label{subsequential domination} Every subsequence of $(t^*_n)$ is $1$-dominated by  $(t^*_n)$ in
$T^*$.

\item\label{upper semi-homogenious} There exists a constant $K$ such that every normalized block
basis in $T^*$ is $K$-dominated by $(t^*_n)$ in $T^*$.

\end{enumerate}

\end{prop}

Statements \ref{equivalence by range} and \ref{subsequential domination} follow from dualizing \cite[Proposition II.4 (b)]{CS} and \cite[Proposition I.9 3.]{CS}  respectively. Statement \ref{upper semi-homogenious} follows directly from the previous two.

\begin{cor}
\label{subsymmetric norm}
The norm $\|\cdot\|$ of $Su(T^*)$ is 1-subsymmetric.
\end{cor}

\begin{proof}
By the definition of $\|\cdot\|$, the basis of $Su(T^*)$ is 1-unconditional and 1-dominated by all of its subsequences. By Proposition \ref{pr2} \ref{subsequential domination}, the reverse domination is also true. 
\end{proof}

\begin{lemma}\label{l3}
Let $(y_n)$ be a normalized block basis in $Su(T^*)$ with $\|y_n\|_\infty \to 0$ as $n\to\infty$. Then $(y_n)$ has a subsequence that is equivalent to a block basis in $T^*$.
\end{lemma}

\begin{proof}
By passing to a subsequence and relabelling we may write
$y_n=\sum_{p_n+1}^{p_{n+1}} a_it_i^*$, where
\begin{equation}
\label{1*}
p_n\max\{|a_i|\colon p_n+1\le|a_i|\le p_{n+1}\}<2^{-(n+1)}, \quad (n\ge 1).
\end{equation}

For $A\subset \N$ let $P_A$ denote the natural basis projection onto $spn\{t_j^*: j\in A\}$. That is, if $x=\sum b_jt_j^*$ then $P_A(x)=\sum_{j\in A} b_jt^*_j$. For the remainder of this proof let also $E_n = \{p_n+1,\ldots,p_{n+1}\}$ and $P_n = P_{E_n}$, $n\geq 1$.

For each $n\ge 1$, note that $1=\|y_n\|=\|\tilde{y}_n\|_{T^*}$, where
$\tilde{y}_n$ is a ``left-spread'' of $P_{A_n}(y_n)$ for some $A_n\subseteq [p_n+1, p_{n+1}]$.
Let $z_n=P_n(\tilde{y}_n)$. Then by \eqref{1*} $\|z_n\|_{T^*}\ge\|\tilde{y}_n\|_{T^*}-1/2=1/2$.

To complete the proof, we show that $(y_n)\subseteq Su(T^*)$ and $(z_n)\subseteq T^*$ are equivalent basic sequences.

For any scalars $(c_n)_{n\ge 1}$ note that the coefficient sequence (with respect to $(t_i^*)$) of $\sum c_nz_n$ is a subsequence of the coefficient sequence of $\sum c_ny_n$. Hence
\[
\left\|\sum c_nz_n\right\|_{T^*}\le \left\|\sum c_nz_n\right\|\le
\left\|\sum c_ny_n\right\|.
\]
Conversely, suppose that $\left\|\sum c_ny_n\right\|=1$. Then $1=\left\|\sum c_nw_n\right\|_{T^*}$, where $w_n$ is a ``left spread'' of $P_{B_n}(y_n)$ for some
$B_n\subseteq [p_n+1,p_{n+1}]$. By \eqref{1*},
$\|w_n-P_n(w_n)\|_{T^*}\le 2^{-(n+1)}$, so
$$
\left\|\sum c_n P_n(w_n)\right\|_{T^*}\ge
\left\|\sum c_nw_n\right\|_{T^*}-\sum_{n=1}^\infty 2^{-(n+1)}=1-\frac12=\frac12.
$$
But since the coefficient sequence of $P_n(w_n)$  is a subsequence of that of  $y_n$
$$
\|P_n(w_n)\|_{T^*}\le\|P_n(w_n)\|\le\|y_n\|=1 \le 2\|z_n\|_{T^*}.
$$
Since $\supp(P_n(w_n))\subseteq [p_n+1,p_{n+1}]$ and
$\supp(z_n)\subseteq [p_n+1,p_{n+1}]$ it follows from Proposition  \ref{pr2} \ref{equivalence by range} that $(\|P_n(w_n)\|_{T^*}^{-1}{P_n(w_n)})$ and
$(\|z_n\|_{T^*}^{-1}{z_n})$ are $K$-equivalent (by perhaps skipping the terms for which $P_n(w_n) = 0$). Thus
$$
\frac12\le\left\|\sum c_nP_n(w_n)\right\|_{T^*}\le
2K\left\|\sum c_nz_n\right\|_{T^*},
$$
i.e.,
$$
\left\|\sum c_nz_n\right\|_{T^*}\ge\frac1{4K}.
$$
Therefore $(y_n)\subseteq Su(T^*)$ and $(z_n)\subseteq T^*$ are equivalent.
\end{proof}

\begin{cor}\label{cor4}
Suppose $(y_n)$ is a
subsymmetric block basis in $Su(T^*)$. Then there exists $\delta>0$ such that, for all $n\ge 1$, $\|y_n\|_\infty>\delta$.
\end{cor}

\begin{proof}
This follows from Lemma \ref{l3} since $T^*$ does not contain any subsymmetric basic sequence.
\end{proof}

\begin{prop}\label{pr5}
Let $(y_n)\subseteq Su(T^*)$ be a normalized block basis. Then $(y_n)$ is $K$-dominated by $(t_n^*)_n$.
\end{prop}

\begin{proof}
Suppose $\left\|\sum c_ny_n\right\|=1$. Therefore $1=\left\|\sum c_nw_n\right\|_{T^*}$ where $w_n$ is a ``left spread'' of $P_{A_n}(y_n)$ for some $A_n\subseteq [p_n+1,p_{n+1}]$. Note that
$$
\|w_n\|_{T^*}\le\|w_n\|\le \|y_n\|=1.
$$
So $\left\|\sum c_nw_n\right\|_{T^*}\le K\left\|\sum c_nt_n^*\right\|_{T^*}\le K\left\|\sum c_nt_n^*\right\|$ (by  Proposition  \ref{pr2} \ref{upper semi-homogenious}), i.e.,
$$
\left\|\sum c_ny_n\right\|\le K\left\|\sum c_nt_n^*\right\|.
$$
\end{proof}

\begin{proof}[Proof of Theorem \ref{th2}]
In every Banach space with a basis, a subsymmetric basic sequence is equivalent to a normalized block basis. Let $y_n=\sum_{p_n+1}^{p_{n+1}} a_it_i^*$, ($n\ge1)$, be a normalized subsymmetric block basis. By Corollary \ref{cor4} there
exists $\delta>0$ such that
$$
\left\| \sum c_ny_n\right\|\ge\delta\|\sum c_nt_n^*\|
$$
and by Proposition \ref{pr5},
$$
\left\|\sum c_ny_n\right\|\le K\left\|\sum c_nt_n^*\right\|.
$$
\end{proof}

\begin{rem}
Note that Theorem \ref{th2} in conjunction with James's  theorem for spaces with an unconditional basis \cite{J} yields that $Su(T^*)$ is reflexive. 
\end{rem}

 \begin{thm} There is a continuum of isomorphically distinct Banach spaces with a subsymmetric basis  (which is not symmetric)  and a unique subsymmetric basic sequence up to equivalence.
\begin{proof}
 Let  $(z_i^*)$ be any subsequence of $(t_i^*)$ and let $Su((z_i^*))$ be the `subsymmetrization' of  $(z_i^*)$. Then $(z_i^*)$ is a subsymmetric basis of 
 $Su((z_i^*))$ with fundamental function $\Phi(n) =\left\|\sum_{i=1}^n
 z_i^*\right\|_{T^*}$.  The proof of Theorem~\ref{th2} easily
 generalizes to show that  $Su((z_i^*))$  has a unique subsymmetric
 basic sequence and that the basis of  $Su((z_i^*))$  is not
 symmetric. It is known
that $(t^*_i)$ has a continuum of non-equivalent subsequences
 \cite{CJT}, \cite{CS} . By  varying the choice of the subsequence $(z_i^*)$ one
 can construct a continuum of such spaces  $Su((z_i^*))$  having
 mutually non-equivalent fundamental functions. Hence the corresponding subsymmetric sequences $(z_i^*)$ in  $Su((z_i^*))$  are not equivalent, which implies that the corresponding spaces  $Su((z_i^*))$  are isomorphically distinct. \end{proof}
 \end{thm}

%%%%%%%%%%%%%%%%%%%%%%%%%%%%%%%%%%%%%%%%%%%%%%%%%%%%%%%%%%%%

\section{$Su(T^*)^*$ has no symmetric basic sequence}
\label{dual empty of symmetry}
In this section we prove that $Su(T^*)^*$ has no symmetric basic sequences. This ought to be contrasted to Corollary \ref{c39}, according to which it contains a continuum of nonequivalent subsymmetric basic sequences.

\begin{lemma}\label{l21}
Let $(e_i)$ be a normalized $1$-subsymmetric basis for a Banach space $X$ not isomorphic to $c_0$. Let $\varepsilon>0$ and let $C>1$. There exists $\delta=\delta(\varepsilon,C)>0$ such that for every $y=\sum_{i=1}^n b_ie_i^*\in X^*$ satisfying $\frac1C\le\|y\|\le C$ and $\max_{i\le i\le n} |b_i|<\delta$ there exists $x=\sum_{i=1}^n a_ie_i\in X$ with $\|x\|\le 1$, $\max_{1\le i\le n} |a_i|<\varepsilon$, and $(x,y)\ge\frac1{2C}$.
\end{lemma}

\begin{proof}
By the Hahn-Banach theorem and 1-unconditionality there exists $z=\sum_{i=1}^n c_ie_i\in X$ such that $(z,y)=\|y\|$ and $\|z\|=1$.
Let $A=\{1\le i\le n \colon |c_i|\ge\varepsilon\}$.
Then $1=\|z\|\ge\varepsilon \Phi(|A|)$, where $(\Phi(n))_{n=1}^\infty$ is the fundamental function of $(e_i)$, i.e. $\Phi(|A|)\le \frac1\varepsilon$.
Since $(e_i)$ is not equivalent to the unit vector basis of $c_0$, $\Phi(m)\to\infty$ as $m\to\infty$. Hence $|A|\le N(\varepsilon)$ for some positive integer $N(\varepsilon)$, so
$$
\left|\sum_{i\in A} b_ic_i\right|\le N(\varepsilon)\max_{1\le i\le n}|b_i|\le N(\varepsilon)\delta.
$$
Setting $\delta=\frac1{2CN(\varepsilon)}$, $B = \{1,\ldots,n\}\setminus A$, and $x= P_Bz$ we have $(x,y)\ge\|y\|-N(\varepsilon)\delta\ge 1/C-2/(2C)=1/(2C)$ and $\|x\|\le\|z\|=1$.
\end{proof}

\begin{lemma}\label{l22}
Suppose $n\ge 1$ and that $x_1,\dots, x_n$ in $Su(T^*)$ satisfy
$\|x_k\|\le 1$, $\|x_k\|_\infty\le\frac1{k2^k}$ and $x_n<x_{n-1}\cdots<x_1$.
Then $\left\|\sum_{k=1}^n x_k\right\|<3$.
\end{lemma}

\begin{proof}
  We prove by induction that
  $\left\|\sum_{k=1}^n x_k\right\|<3-2^{-n}$. Suppose the result holds for
  a natural number
  $n$ and all $x_1, x_2,\cdots,x_n$ . Let us prove the result for $n+1$. So assume
$x_{n+1}<x_n<\cdots<x_1$.
Let $\sum_{k=1}^{n+1}x_k=\sum a_ie_i$.
Then
\[\left\|\sum_{k=1}^{n+1} x_k\right\|=\left\|\sum a_{m_i}e_i\right\|_{T^*}\text{ for some }m_1<m_2<\cdots.\]
We consider two cases.

First, if $m_{n+1}\le \max\supp(x_{n+1})$, then
$$
\left\|\sum_{k=1}^{n+1}x_k\right\|\le (n+1)\|x_{n+1}\|_{\infty}+
\left\|\bar{x}_{n+1}+\sum_{k=1}^n\bar{x}_k\right\|_{T^*}
$$
(where $n+2\le \bar{x}_{n+1}<\bar{x}_n<\cdots<\bar{x}_1$ and $\|\bar{x}_k\|_{T^*}\le 1$)
$$
\le (n+1)\frac1{2^{n+1}(n+1)}+2
$$
(since $n+1<\bar{x}_{n+1}<\bar{x}_n<\cdots<\bar{x}_1$
and $\|\bar{x}_k\|_{T^*}\le 1$ for $1\le k\le n+1$)
$$
=2+\frac1{2^{n+1}}<3-\frac1{2^{n+1}}.
$$

In the second case, $m_{n+1}>\supp(x_{n+1})$. Then
$$
\left\|\sum_{k=1}^{n+1}x_k\right\|\le n\|x_{n+1}\|_\infty+\left\|\sum_{k=1}^n x_k\right\|
\le\frac{n}{(n+1)2^{n+1}}+3-2^{-n}
$$
(by the induction hypothesis)
$$
<3-2^{-(n+1)}.
$$
\end{proof}

\begin{thm}\label{th21}
$Su(T^*)^*$ does not contain any symmetric basic sequence.
\end{thm}

\begin{proof}
Suppose, to derive a contradiction, that $(z_i)$ is a semi-normalized symmetric basic sequence in $Su(T^*)^*$. We may assume that $(z_i)$ is a block basis with respect to $(e_i^*)$. We may also assume that $\|z_i\|_\infty\le 1$ for each $i$. Let $C$ be the symmetry constant of $(z_i)$. For each $n\ge 1$, let $\delta_n=\delta_n(\varepsilon_n, C)$, where $\varepsilon_n=\frac{1}{n2^n}$. Since $(z_i)$ is not equivalent to the unit vector basis of $c_0$ there exists a normalized block basis (of $(z_i)$)
$y_k=\sum_{i=1}^{p_k-p_{k-1}} c_i^k z_{p_{k-1}+i}$ such that $\max_{1\le i\le p_k-p_{k-1}}|c_i^k|<\delta_k$ for each $k\ge 1$.
(Here $0=p_0<p_1<p_2<\cdots$.)

For each $n\ge 1$ and $1\le k\le n$, let $y_k^n=\sum_{i=1}^{p_k-p_{k-1}} c_i^k z_{i+p_n-p_k}$. Note that $y_n^n<y_{n-1}^n<\cdots<y_1^n$, and that $y_k^n$ is a translate of $y_k$ relative to the basis $(z_i)$.

Since $\|y_k\|=1$, the symmetry of $(z_i)$ gives $\frac1C\le\|y_k^n\|\le C$.
Note also that $\|y_k^n\|_\infty<\delta_k$ since $\|z_i\|_\infty\le 1$.

By Lemma \ref{l21} there exist $x_k^n\in Su(T^*)$ with $\|x_k^n\|_\infty\le\varepsilon_k$, $\|x_k^n\|\le 1$, and $(x_k^n,y_k^n)\ge\frac1{2C}$.

Clearly, we may also assume that $\supp(x_k^n)\subseteq\supp(y_k^n)$ by 1-uncondi\-tion\-al\-ity of $(e_i)$.

By Lemma \ref{l22}, since
$\|x_k^n\|_\infty\le\varepsilon_k=\frac{1}{k2^k}$, we obtain
$\left\|\sum_{k=1}^n x_k^n\right\|<3$.

Hence, for all scalars $a_1,\dots,a_n$
\begin{eqnarray*}
\left\|\sum_{k=1}^n a_ky_k^n\right\|&\ge&\frac13
\left(\sum_{k=1}^n\sgn(a_k)x_k^n, \sum_{k=1}^n a_ky_k^n\right)\\
&=&\frac13\sum_{k=1}^n |a_k|(x_k^n,y_k^n)\ge\frac1{6C}\sum_{k=1}^n |a_k|.
\end{eqnarray*}
On the othed hand,
$$
\left\|\sum_{k=1}^n a_ky_k^n\right\|\le
\sum_{k=1}^n|a_n|\|y_k^n\|\le C\sum_{k=1}^n|a_k|.
$$
By the symmetry of $(z_i)$, it follows that for each $n$
$$
\frac1{6C^2}\sum_{k=1}^n|a_k|\le\left\|\sum_{k=1}^na_ky_k\right\|\le C^2\sum_{k=1}^n|a_k|.
$$
So $(y_k)_{k=1}^\infty$ is equivalent to the unit vector basis of $\ell_1$. This is impossible since $Su(T^*)^*$ is reflexive.
\end{proof}

By duality we also obtain the following.

\begin{cor}\label{c21}
$Su(T^*)$ does not have any quotient space with a symmetric basis.
\end{cor}

\section{Subsymmetric bases equivalent to the unit vector bases of $c_0$ or $\ell_p$}

In this section we extend to subsymmetric bases the main result of \cite{A1}. This result also ought to be compared to \cite[Proposition 7]{FPR}.

\begin{thm}
\label{l_p characterization}
Let $(e_i)$ be a subsymmetric basis of a Banach space $X$. Suppose that $(e_i)$ dominates every subsymmetric block basis  with respect to $(e_i)$ and that $(e_i^*)$ dominates every subsymmetric  block basis 
with respect to $(e_i^*)$. Then $(e_i)$ is equivalent to the unit vector basis of $c_0$ or $\ell_p$ for some $1 \le p < \infty$.
\end{thm}

The proof requires some notation and lemmata.

~
\begin{enumerate}[leftmargin=19pt]

\item Let $X$ be a space with a 1-subsymmetric normalized basis $(e_i)_{i=1}^\infty$.

\item Let $X(\omega^2)$ be the closure of  the vector space (with basis $\{e_\alpha \colon \alpha < \omega^2\}$) of  finitely supported vectors $\sum_{i=1}^n a_i e_{\alpha_i}$,
where $0 \le \alpha_1<\alpha_2<\dots<\alpha_n < \omega^2$,  equipped with the norm

$$
\left\|\sum_{i=1}^n  a_i e_{\alpha_i}\right\|= \left\|\sum_{i=1}^n  a_i e_i\right\|.
$$
Each $x$ in the space $X(\omega^2)$ has a unique representation as a formal series $\sum_{\alpha < \omega^2} a_\alpha e_\alpha$.

\item For $0\ne x=\sum_{j=1}^\infty a_je_j \in X$,  let $x^i=\sum_{j=1}^\infty a_j e_{\omega(i-1)+j}\in X(\omega^2)$ for $i\ge 1$.

\end{enumerate}
\begin{lemma}\label{l35}
$(x^i)_{i=1}^\infty$ is $1$-subsymmetric in $X(\omega^2)$ and is equivalent to a subsymmetric sequence $(y_i)_{i=1}^\infty$ in $X$.
\end{lemma}

\begin{proof} The first assertion is clear from the definition of the norm in $X(\omega^2)$ and the fact that  $(e_i)$ is $1$-subsymmetric.
Choose finitely supported vectors  $\tilde{y}_i$ in $X(\omega^2)$ with $\operatorname{supp} (\tilde{y}_i) \subseteq \operatorname{supp} (x^i )$ and
$\|x^i-\tilde{y_i}\|<\varepsilon_i\downarrow0$. Now choose
$y_1<y_2<\cdots$ in $X$ such that $y_i$ has the same distribution as
$\tilde{y}_i$. Then $(x_i)$ and $(y_i)$ are equivalent basic sequences  provided $\varepsilon_i\downarrow0$ sufficiently fast.
\end{proof}
\begin{lemma}
Let  $(e_i)$ be a normalized  $1$-subsymmetric basis for a Banach space $X$. % that dominates every subsymmetric block basis of itself. 
Suppose that for each $x \in X$ there exists a constant $K(x)>0$ such that  $(x^i) \subset X(\omega^2)$ satisfies,  for each $y = \sum a_i e_i \in X$, 
 $$\left\| \sum a_i x^i\right\| \le K(x) \|y\|.$$
Then there exists $K>0$ such that for all $x\in X$, $K(x)\leq K\|x\|$.
\end{lemma}
\begin{proof} A standard gliding hump argument shows that, for each
$y = \sum a_i e_i \in X$,  there exist $C(y) > 0$ such that 
$$\left\| \sum a_i x^i\right\| \le C(y) \|x\|.$$
Hence the linear  mapping $T_y \colon X \rightarrow X(\omega^2)$ given by
$T_y(x) = \sum a_i x^i$ is continuous. Moreover, for each $x \in X$, 
$$\sup \{ \|T_y(x) \|\colon \|y\| \le 1 \} \le K(x).$$
The desired conclusion now follows from the Uniform Boundedness Principle.
\end{proof}

Note that $(e^*_i)$ is $1$-subsymmetric in $X^*$. Let  $Y \subseteq X^*$ be the closed linear span of $(e^*_i)$.
Note that $Y(\omega^2)$  is isometrically isomorphic to a subspace of $X(\omega^2)^*$ with the  duality pairing of $Y(\omega^2) \times X(\omega^2)$
given by 
$$\langle \sum b_\alpha e_\alpha, \sum a_\alpha e_\alpha \rangle= \sum b_\alpha a_\alpha$$

\begin{lemma}
Let  $(e_i)$ be a normalized  $1$-subsymmetric basis for a Banach space $X$. 
%and assume that $(e_i^*)$ dominates every subsymmetric block basis of itself. 
Suppose that for each $f \in Y$ there exists a constant $K(f)>0$ such that  $(f^i) \subset Y(\omega^2)$ satisfies,  for each $g = \sum b_i e^*_i \in Y$, 
 $$\left\| \sum b_i f^i\right\| \le K(f) \|g\|.$$
Then there exists $C>0$ such that  for all $x\in X$ and $y = \sum a_i e_i \in X$,
$$\left\| \sum a_i x^i\right\| \ge  C\|x\| \|y\|.$$
\end{lemma} 
\begin{proof} By the previous lemma, there exists $K>0$ such that, for all $f\in Y$ and $g = \sum b_i e_i \in Y$,
 $$\left\| \sum b_i f^i\right\| \le K \|f\|  \|g\|.$$
Suppose that $x\in X$, $y = \sum a_i e_i \in X$ and that $\|x\|=\|y\|=1$. Select $f \in Y$, with $\|f\|=1$ and $f(y)>1/2$, and select
$g = \sum b_i e^*_i \in Y$ with $\|g\|=1$ and $g(f) = \sum |a_i| |b_i| > 1/2$. Then $\left\|\sum b_i f ^i\right\| \le K$ and 
$$\langle \sum b_i f^i, \sum a_i x^i \rangle  \ge \frac{1}{2} \sum |a_i| |b_i| \ge \frac{1}{4}.$$
Hence $\left\|\sum a_i x^i\right\| \ge 1/4K$. The result for $C= 1/4K$ follows by homogeneity.
\end{proof}

\begin{proof}[Proof of Theorem \ref{l_p characterization}] For each  $0 \ne x \in X$,  $(x^i)$ is equivalent to a subsymmetric block basis with respect to $(e_i)$. Similiarly, for each $0 \ne f \in Y$, $(f^i)$ is equivalent to a subsymmetric block basis with respect to $(e_i^*)$. Hence the hypotheses of the previous two lemmas are satisfied. It follows that there exists $K< \infty$ such that for all $x \in X$ and $y = \sum a_ie_i \in X$,
$$\frac{1}{K} \|x\| \|y\| \le \left\|\sum a_i x^i\right\| \le K \|x\| \|y\|.$$ For $m,n \in \mathbb{N}$, setting $x =  \sum_{i=1}^m e_i$ and $y = \sum_{i=1}^n e_i$
yields $$\frac{1}{K} \Phi(m) \Phi(n) \le \Phi(mn) \le K \Phi(m) \Phi(n),$$
where $(\Phi(n))$ is the fundamental function of $(e_i)$.
The proof is now concluded as in \cite{A1} (or \cite{LT2}).  It
suffices to observe that the argument presented in \cite{A1} for a
symmetric basis depends only on  the subsymmetry of the basis. (Alternatively, the proof can be concluded with  an argument based on Krivine's theorem \cite{K} as in  \cite[Theorem~1]{FPR}.)
\end{proof}
\vskip12pt

\begin{rem}
  As in the result of Altshuler, it was sufficient to work only with
  subsymmetric block bases which are equivalent to subsymmetric
  sequences in $X(\omega^2)$ and $Y(\omega^2)$ generated by a vector.
 \end{rem}
\section{Spaces with a continuum of subsymmetric sequences}In this section we give criteria for spaces with a subsymmetric basis to contain a continuum of non-equivalent subsymmetric sequences. They can be applied, e.g., to $Su(T^*)^*$ and Schlumprecht space.
 \begin{thm} \label{thm: uncountable}
 Let $(e_i)$ be a $1$-subsymmetric basis for a Banach space $X$.  Suppose that there exists a unit vector $x \in X$ such that
$(e_i)$ does not dominate  $(x^i) \subseteq X(\omega^2)$. Then $(e_i)$ admits a continuum of non-equivalent subsymmetric block bases.
\end{thm} \begin{proof} We shall construct inductively the following: \begin{enumerate} \item  a  normalized block basis  $(x_n)$ of $(e_i)$;
\item an increasing sequence $(N(n))_{n=1}^\infty$ of positive integers;
\item a sequence $y_n = \sum_{i=1}^{N(n) } b_{n,i} e_i$ ($n \ge 1$)   of unit vectors in $X$;
\item a postive sequence $\delta_n \rightarrow 0$.\end{enumerate}
To start the induction, set  $x_1 =y_1=  e_1$, $N(1)=1$,  and $\delta_1 = 1/2$. Suppose that $N>1$ and that $x_k, y_k, N(k)$ and $\delta_k$ have been defined for $1 \le k \le n-1$. 

Let $y_{n-1} = x_1 + x_2 +\dots+x_{n-1}$. Since $y_{n-1}$ has finite support  with respect to $(e_i)$, it follows that there exists $K_n>0$ such that
\begin{equation} \label{eq: K_n}
 \|\sum a_i y_{n-1}^i\|_{X(\omega^2)}  \le K_n \|\sum a_i e_i\| \end{equation} for all coefficient sequences  $(a_i)$. Let \begin{equation}\label{eq: delta_n}
  \delta_n =\frac{ 2^{-n}}{N(n-1)}. \end{equation}
Since $(e_i)$ does not dominate $(x^i) \subset X(\omega^2)$, there exist a positive integer  $N(n)>N(n-1)$ and a unit vector  $y_n = \sum_{i=1}^{N(n) } b_{n,i} e_i$ such that
$$\|\sum_{i=1}^{N(n)} b_{n,i} x^i\|_{X(\omega^2)} >  \frac{n K_n}{\delta_n}.$$
Let  $z_n\in X$ be a finitely supported (with respect to $(e_i)$) approximation to $X$  and let  $x_n>x_{n-1}$ be a right shift of $z_n/\|z_n\|$. Then
$\|x_n\|=1$ and, provided  $\|x - z_n\|$ is  sufficiently 
small, \begin{equation} \label{eq: lowerbound}
\|\sum_{i=1}^{N(n)} b_{n,i} x_n^i\|_{X(\omega^2)} >  \frac{n K_n}{\delta_n}. \end{equation}This completes the inductive step.

Now suppose that $\varepsilon=(\varepsilon_n)_{n=1}^\infty\in \{0,1\}^{\mathbb{N}}$
and that $\varepsilon_1=1$.

Let $x_\varepsilon=\sum_{n=1}^\infty \varepsilon_n\delta_nx_n$. Then
$$
\frac{1}{2} =\varepsilon_1 \delta_1 \le  \|x_\varepsilon\|\le\sum_{n=1}^\infty\delta_n\le\sum_{n=1}^\infty 2^{-n}=1.
$$

Consider the subsymmetric sequence $(x^i_\varepsilon)_{i=1}^\infty$ in $X(\omega^2)$.

($\bullet$) Suppose $\varepsilon_n=1$. Then by \eqref{eq: lowerbound}
\begin{equation} \label{eq: nK_n}
\| \sum_{i=1}^{N(n)} b_{n,i} x^i_\varepsilon\|_{X(\omega^2)} \ge \delta_n \|\sum_{i=1}^{N(n)} b_{n,i} x_n^i\|_{X(\omega^2)}> nK_n
\end{equation}

($\bullet\bullet$) Suppose $\varepsilon_n=0$. Note that $|b_{n,i}|\le 1$ for all $n$ and $i$ and that $K_n \ge 1$ by definition. So by \eqref{eq: K_n} and \eqref{eq: delta_n}
\begin{equation} \label{eq: upperbound}  \begin{split}
\| \sum_{i=1}^{N(n)} b_{n,i} x^i_\varepsilon\|_{X(\omega^2)} \le&\| \sum_{i=1}^{N(n)} b_{n,i} y^i_{n-1}\|_{X(\omega^2)} 
+\sum_{i=n+1}^\infty \delta_i N(n)\\&\le  K_n
+\sum_{i=n+1}^\infty \frac{1}{2^{i}} \\
&\le 2K_n \end{split}
\end{equation}

To complete the proof, pick a subset $(a_r)_{r\in\mathbb{R}}$ of $\{0,1\}^\N$ such that any $r<s$, $\alpha_r = (\alpha_n),\alpha_s = (\beta_n)$ satisfy the following property: $\alpha_n\le\beta_n$ for all $n$ and $\alpha_n<\beta_n$ for infinitely many $n$ (and $\alpha_1=1$). This can be achieved by enumerating $\mathbb{Q} = (q_n)_n$ and putting for each $r\in (q_1, \infty)$, $A_r = \{n:q_n<r\}$ and taking its characteristic $\alpha_r$. In conclusion, for $r < s$, $\alpha_r = (\alpha_n),\alpha_s = (\beta_n)$, and $n\in A_s\setminus A_r$, \eqref{eq: nK_n} and \eqref{eq: upperbound} imply  $$
\| \sum_{i=1}^{N(n)} b_{n,i} x^i_{\alpha_s}\|_{X(\omega^2)} \ge \frac{n}{2}\| \sum_{i=1}^{N(n)} b_{n,i} x^i_{\alpha_r}\|_{X(\omega^2)}$$

In particular, $(x^i_{\alpha_r})_i$ and $(x^i_{\alpha_s})_i$ are non-equivalent.

\end{proof}
\begin{remark} \label{rem: dominationorder} Note that if $r < s$ then  $(x^i_{\alpha_r})_i \prec (x^i_{\alpha_s})_i$ in the domination ordering $\preceq$ of basic sequences. So the continuum of subsymmetric block bases constructed in Theorem~\ref{l_p characterization} (and in subsequent results)
is a chain that is
is order-isomorphic to $(\mathbb{R}, \le)$ in the domination ordering.
\end{remark}
The next result is a significant  strengthening of Theorem~\ref{l_p characterization}.
 \begin{thm} \label{thm: l_p char2} Let $(e_i)$ be a subsymmetric basis which is not equivalent to the unit vector basis of $\ell_p$ or $c_0$. Then either $(e_i)$ or $(e_i^*)$ admits a continuum of non-equivalent subsymmetric block bases.
\end{thm} \begin{proof} By Theorem~\ref{l_p characterization}  either $(e_i)$ or $(e_i^*)$ satisfies the hypothesis and hence the conclusion  of Theorem~\ref{thm: uncountable}.   Suppose that $(e_i)$ does. By Lemma~\ref{l35} each $(x_\varepsilon^i)$ is equivalent to a subsymmetric block basis of $(e_i)$. Hence $(e_i)$ admits a continuum of non-equivalent subsymmetric block bases. The same reasoning applies to   $(e_i^*)$.
\end{proof}
The next result is a criterion for the existence of a continuum of subsymmetic sequences with  non-equivalent fundamental functions.
\begin{thm}\label{ell1 by vector}
Suppose that $(e_i)_{i=1}^\infty$ is a 1-subsymmetric normalized basis which is not equivalent to the unit vector basis  of $\ell_1$ and that, for each $n\ge 1$, there exist unit vectors $z_1<z_2<\cdots<z_n$ such that $z_i$ is a right shift of $z_1$ $(1\le i\le n)$ and $\left\|\sum_{i=1}^n z_i\right\|\ge\frac n2$. Then $X$ contains a continuum of subsymmetric sequences with non-equivalent fundamental functions.
\end{thm}\begin{proof}  For $n\ge1$, let $\Phi(n)=\| \sum_{i=1}^n e_i\|$ be the fundamental function of $(e_i)$. Since $(e_i)$ is subsymmetric and not equivalent to the unit vector basis of $\ell_1$, it follows thar $\Phi(n)/n \rightarrow 0$ as $n \rightarrow \infty$. Hence for each $n \ge 1$ there exist 
$N(n) \in \mathbb{N}$  and  unit vectors $z^n_1<z^n_2<\dots z^n_{N(n)}$ such that 
$$\| \sum_{i=1}^{N(n)} z^n_i\| \ge \frac{N(n)}{2} \ge  n2^n \Phi(N(n)).$$
By subsymmetry, we may assume that $z^1_1 < z^2_1 <\dots$. Let $x = \sum_{n=1}^\infty 2^{-n} z^n_1 \in X$.  Then, for each $n \in \mathbb{N}$, 
$$ \| \sum_{i=1}^{N(n)} x^i\|_{X(\omega^2)}\ge 2^{-n} \|\sum_{i=1}^{N(n)} z^n_i\| \ge n\Phi(N(n)).$$
So $(e_i)$ does not dominate $(x^i)$. By Theorem~\ref{thm: uncountable} $(e_i)$ admits a continuum of non-equivalent symmetric block bases.   However, straightforward modifications to the proof of Theorem~\ref{thm: uncountable} yield a stronger result in this case, namely  a continuum of subsymmetric block bases with non-equivalent fundamental functions. 
\end{proof}

\begin{lemma}\label{l32}
For every $n\ge 1$ there exist unit vectors $x_1<x_2<\cdots<x_n$ in $Su(T^*)$, where $x_i$ is a right shift of $x_1$ $(1\le i\le n)$ such that for all $a_1,\dots,a_n$
$$
\max|a_i|\le
\left\|\sum_{i=1}^n a_ix_i\right\|\le 3\max|a_i|.
$$
\end{lemma}

\begin{proof}
For $n\in\N$ pick a finitely supported unit vector $x_1\in X$ with $\|x\|_\infty\leq (n2^n)^{-1}$. This is possible because $(e_i)$ is not equivalent to the unit vector basis of $c_0$. Let $x_1<x_2<\cdots<x_n$ be copies of $x_1$. By Lemma \ref{l22}, $\|\sum_{i=1}^nx_i\|\leq 3$. The result follows from 1-unconditionality of the basis.
\end{proof}

\begin{lemma}\label{l33}\label{r34}
For each $\varepsilon>0$ and $n\ge 1$, $Su(T^*)^*$ contains blocks $x_1<x_2<\cdots<x_n$, where each $x_i$ is a right shift of $x_1$, and for all $a_1,\dots,a_n$,
$$
\frac1{3}\sum_{i=1}^n |a_i|\le\left\|\sum_{i=1}^n a_ix_i\right\|\le
\sum_{i=1}^n |a_i|.
$$
In particular,
$$
\left\|\sum_{i=1}^n x_i\right\|\ge\frac n3.
$$
\end{lemma}

\begin{proof}
This follows from Lemma \ref{l33} by duality.
\end{proof}

\begin{cor}\label{c39}
$Su(T^*)^*$ contains a continuum of  subsymmetric basic sequences with non-equivalent fundamental functions. (Hence, by duality, $Su(T^*)$ has a continuum of  subsymmetric  quotient spaces with non-equivalent fundamental functions.)
\end{cor}

\begin{proof}
The hypothesis of Theorem \ref{ell1 by vector} is satisfied by Lemma \ref{r34}.
\end{proof}

\begin{cor}\label{c312}\cite{AS}
 Schlumprecht space $S$ contains a continuum  of
 subsymmetric basic sequences with non-equivalent fundamental functions.
\end{cor}

\begin{proof}
 The space $S$ satisfies the hypothesis of Corollary \ref{ell1 by vector}, see, e.g.,
  \cite[Lemma 3]{GM}.
\end{proof}
Next we present  another  criterion for a space with a subsymmetric basis to contain a continuum of subsymmetric 
basic sequences with non-equivalent fundamental functions. 
Let us say that the fundamental function $(\Phi(n))$ of $(e_i)$ is \textit{submultiplicative} if there exists $K>0$ such that $\phi(nm) \le K \Phi(n)\Phi(m)$ for all $m,n \in \mathbb{N}$.
\begin{thm}\label{th36}
Suppose that $(e_i)_{i=1}^\infty$ is a 1-subsymmetric  normalized basis whose fundamental function is not submultiplicative. Then $(e_i)$ admits a continuum of subsymmetric block bases with non-equivalent fundamental functions.
\end{thm} \begin{proof} Let $n \in \mathbb{N}$.  Since $(\Phi(r))_{r=1}^\infty$ is not submultiplicative, there exist $r_n,s_n \in \mathbb{N}$ such that
$$\Phi(r_ns_n) >n2^n\Phi(r_n)\Phi(s_n).$$ Hence there exists a normalized block basis $(x_n)$ such that 
$x_n$ is a right shift of $\Phi(r_n)^{-1}\sum_{i=1}^{r(n)} e_i$ and  $\|\sum_{i=1}^{s(n)} x_n^i\|_{X(\omega^2)} > n 2^n \Phi(s_n)$.
Let $x = \sum_{n=1}^\infty 2^{-n} x_n \in X$. Then, for each $n \in \mathbb{N}$,
$$\|\sum_{i=1}^{s(n)} x^i\|_{X(\omega^2)} \ge 2^{-n} \|\sum_{i=1}^{r(n)}x_n^i\|_{X(\omega^2)} \ge n \Phi(s(n)).$$
So $(e_i)$ does not dominate $(x^i)$. Straightforward modifications to the proof of Theorem~\ref{thm: uncountable}  yield a continuum of subsymmetric block bases with non-equivalent fundamental functions.
\end{proof}

\begin{remark} All of the results of this section  remain valid if `subsymmetric' is replaced by `symmetric'.  \end{remark}

\begin{cor}Let $(e_i)$ be the symmetric basis constructed by  Altshuler \cite{A} which is equivalent to all its  symmetric block bases.
Then $(e^*_i)$  admits a continuum of symmetric  block bases with non-equivalent fundamental functions.
\end{cor}

\begin{proof} Since $(e_i)$ is equivalent to its symmetric block bases, it follows from the symmetric version of Theorem~\ref{th36} that its fundamental function 
$(\Phi(n))$ is submultiplicative.
  Let $(\Phi^*(n))$ be the fundamental function of $(e_i^*)$. Suppose, to derive a contradiction,  that  $(\Phi^*(n))$ is  submultiplicative. Since $(e_i)$ is symmetric,
$\Phi^*(n) \asymp n/\Phi(n)$ (see \cite[Prop.~ 3.a.6]{LT1}). Hence there exists $K>0$ such that
$$\frac{1}{K} \Phi(m) \Phi(n) \le \Phi(mn) \le K  \Phi(m) \Phi(n)$$
for all $m,n \in \mathbb{N}$. It follows that $\Phi(n) \asymp n^{1/p}$ for some $p \in [1,\infty]$ (see \cite[Theorem 2.a.9]{LT1}).

However,  using the definition of the norm given in \cite{A1}, $$\Phi(n) \le \sum_{j=1}^n 1/j \le 1 + \ln(n).$$  Moreover, $\Phi(n) \rightarrow \infty$
since $(e_i)$ is not equivalent to the unit vector basis of $c_0$. So $(\Phi(n))$ is not equivalent to $(n^{1/p})$ for any $p \in [1,\infty]$, which  contradicts the assumption. Hence $(\Phi^*(n))$ is not submultiplicative and the result follows from  the  symmetric version of  Theorem~\ref{th36}.
\end{proof}

Next we present  an application of the symmetric version of Theorem~\ref{thm: uncountable} to the classical Lorentz sequence spaces $\ell_{p,q}$.
Let $1 \le p,q < \infty$. Recall that the Lorentz sequence space $\ell_{p,q}$ is the closure of $c_{00}$ under the quasi-norm
$$\|\sum a_i e_i\|_{p,q} = (\sum( a_i^*)^q i^{q/p-1})^{1/q},$$ where $(a_i^*)$ is the non-increasing rearrangement of $(|a_i|)$.
For $1 \le q \le p$, $\|\cdot\|_{p,q}$ is a norm.  However, for $1< p<q<\infty$,  $\|\cdot\|_{p,q}$ does not satisfy the triangle inequality but is nevertheless  equivalent, under the natural duality,  to the dual norm $\|\cdot\|_{p^\prime, q^\prime}^*$,
where $1/p + 1/p^\prime = 1/q + 1/q^\prime =1$.

Part (b) of the following corollary appears to be new to the best of our knowledge.

\begin{cor} (a)   For $1 \le q < p$, $\ell_{p,q}$ contains exactly two non-equivalent symmetric basic sequences, viz., the unit vector bases of $\ell_{p,q}$ and of $\ell_q$.
\newline
(b) For $1< p<q<\infty$, the unit vector basis of  $\ell_{p,q}$ admis a continuum of non-equivalent symmetric   block bases. 
\end{cor}  \begin{proof}  (a) This follows from \cite[Theorem~6]{ACL} since $\ell_{p,q} = d(a,q)$ with submultiplicative weight $a = (i^{q/p-1})$. \newline
(b) This follows by combining (a) with the symmetric version of Theorem 19 since $\ell_{p,q} =  \ell_{p^\prime, q^\prime}^*$ with an equivalent norm.
\end{proof}
\section{Applications to spreading models}
In  this section we apply our results  to the setting of spreading models. We thusly obtain a criterion for a Banach space to admit a continuum of pairwise non-equivalent spreading models.

A Schauder basic sequence $(x_i)$ in a Banach space $X$ is said to generating a sequence $(e_i)$ in a Banach space $E$ as spreading model if for any $a_1,\ldots,a_n$ in $\mathbb{R}$,
\[\lim_{i_1\to\infty}\cdots\lim_{i_n\to\infty}\big\|\sum_{k=1}^na_k x_{i_k}\big\|_X = \big\|\sum_{k=1}^na_k e_{k}\big\|_E.\]
Up to passing to a subsequence of $(x_i)$, the above is equivalent to saying that for a pre-chosen null sequence of positive real numbers $(\delta_n)$, for every $a_1,\ldots,a_n$ in $[-1,1]$, and for every $n\leq i_1<\cdots<i_n$,
\begin{equation}
\label{sm quantitative}
\Big|\big\|\sum_{k=1}^na_k x_{i_k}\big\|_X - \big\|\sum_{k=1}^na_k e_{k}\big\|\Big|_E\leq \delta_n.
\end{equation}
The spreading model $(e_i)$ is always 1-spreading. If we additionally assume that $(x_i)_i$ is weakly null then it is $1$-suppression unconditional.

\begin{prop}
\label{vector generated sm}
Let $(x_i)$ be a Schauder basic sequence in a Banach space $X$ generating a subsymmetric spreading model $(e_i)$ and let $u = \sum_{i=1}^\infty c_ie_i$ be a non-zero vector in the closed linear span of $(e_i)$. Then, there exists a block sequence of $(x_i)$ that generates $(u^i)$ as a spreading model.
\end{prop}

\begin{proof}
We may assume, without loss of generality, that for each $i\in\N$, $c_i\in[-1,1]$ and that $(x_i)$ satisfies \eqref{sm quantitative} for $\delta_n = 1/2^{n}$. Pick $p_1<p_2<\cdots$ so that for all $n\in\mathbb{N}$,
\[\Big\|\sum_{i=p_n+1}^{p_{n+1}}c_ie_i\Big\| \leq 1/2^{n}\]
and define
$u_n = \sum_{i=1}^{p_n}c_ie_i$. Then, for each $n\in\mathbb{N}$ and $a_1,\ldots,a_n\in[-1,1]$ we have
\[\Big|\big\|\sum_{i=1}^na_iu_n^i\big\| - \big\|\sum_{i=1}^na_iu^i\big\|\Big| \leq 2n/2^{n}.\]
We choose a block sequence $(y_n)$ so that for each $n\in\N$,
\[y_n = \sum_{i=1}^{p_n}c_ix_{j^n_i}\]
with $np_n \leq j_1^n<\cdots<j_{p_n}^n$. That is, $y_n$ has the same distribution as $u_n$ and its support starts after a sufficiently large number. 

We claim that $(y_k)$ generates $(u^i)$ as a spreading model. First observe that for any $k\leq n$, if we put $y_n^{(k)} = \sum_{i=1}^{p_k}c_ix_{j^n_i}$ then
\[\big\|y_n - y_n^{(k)}\big\| \leq \|u_n - u_k\| + 1/2^{n} \leq 2/2^k + 1/2^n\leq 3/2^k\]
This is because $y_n - y_n^{(k)}$ has the same distribution as $u_n - u_{k}$ and it starts after $p_n$. Similarly, if we fix $a_1,\ldots,a_n\in[-1,1]$ and $n\leq k_1<\cdots<k_n$ then
\[\Big|\big\|\sum_{i=1}^n a_i y_{k_i}^{(n)}\big\| - \big\|\sum_{i=1}^n a_i u_n^i\big\|\Big| \leq 1/2^{n}.\]
This is because $\sum_{i=1}^n a_i y_{k_i}^{(n)}$ has the same distribution as the vector $\sum_{i=1}^n a_i u_n^i$, all its coefficients (relative to the basis $(x_i)$) are in $[-1,1]$, and its support has at most $np_n$ members that start after $np_n$.

All that remains is to apply the triangle inequality.
\begin{align*}
\Big|\big\|\sum_{i=1}^n a_i y_{k_i}\big\| - \big\|\sum_{i=1}^n a_i u^i\big\| \Big| &\leq \Big|\big\|\sum_{i=1}^n a_i y_{k_i}^{(n)}\big\| - \big\|\sum_{i=1}^n a_i u_n^i\big\|\Big|\\
&+ \sum_{i=1}^n\|y_{k_i} - y_{k_i}^{(n)}\| + \Big|\big\|\sum_{i=1}^ka_iu_n^i\big\| - \big\|\sum_{i=1}^na_iu^i\big\|\Big|\\
&\leq 1/2^n+ 3/2^{n} + 2n/2^n \leq 5n/2^n.
\end{align*}
\end{proof}

\begin{rem}
The assumption that $(e_i)$ is subsymmetric is not essential. The same proof would work if $(e_i)$ were merely a spreading Schauder basic sequence. The space $E(\omega^2)$ can in this case be defined just as well.
\end{rem}

\begin{prop}
\label{sm unc}
Let $X$ be a Banach space that admits a subsymmetric spreading model $(e_i)$ that satisfies one of the following properties.
\begin{enumerate}[label=(\roman*)]

\item\label{sm unc1} There exists a unit vector $z$ in the closed linear span of $(e_i)$ such that $(e_i)$ does not dominate $(z^i)$.

\item\label{sm unc2} For every $n\in\N$ there exists a unit vector $z$ in the closed linear span of $(e_i)$ such that $\|\sum_{i=1}^nz^i\| \geq n/2$.

\item\label{sm unc3} The fundamental function of $(e_i)$ is not submultiplicative.

\end{enumerate}
Then $X$ admits a continuum of mutually non-equivalent spreading models. If \ref{sm unc2} or \ref{sm unc3} holds then $X$ admits a continuum of spreading models with mutually non-equivalent fundamental functions.
\end{prop}

\begin{proof}
In either of these three cases, by applying Theorem \ref{thm: uncountable}, Theorem \ref{ell1 by vector}, or Theorem \ref{th36}, the space $E = \overline{\spn(e_i)}$ contains a continuum of non-zero vectors $(u_\alpha)_{\alpha<\mathfrak{c}}$ so that the subsymmetric sequences $(u^i_\alpha)$ are mutually non-equivalent. By Proposition  \ref{vector generated sm} the result follows.
\end{proof}

For a Banach space $X$ with a basis $(e_i)$ denote by $SP_{w}(X)$ the collection of Schauder basic sequences that are generated as a spreading model by some normalized weakly null sequence in $X$. As it was noted earlier, all such sequences are 1-suppression unconditional. Recall that a Schauder basic sequence $(e_i)$ in a Banach space $X$ is called $C$-Schreier unconditional, for a constant $C>0$, if for every Schreier set $F$ (i.e., $F \subset \mathbb{N}$ with $\min(F) \geq |F|$) and coefficients $(a_i)$ we have
\[\Big\|\sum_{i\in F}a_ie_i\Big\| \leq C\Big\|\sum_{i=1}^\infty a_ie_i\Big\|.\]
In other words, the canonical basis projection from $\overline{\spn(e_i)}$ to $\spn(e_i)_{i\in F}$ has norm at most $C$. Obviously, this property is then also enjoyed by $(e_i^*)$.  Additionally, if a subsequence of $(e_i)$ generates a spreading model $(z_i)$ then a subsequence of $(e_i^*)_i$ generates a spreading model that is $C$-equivalent to $(z_i^*)$. This is because for every Schreier set $F$, $(e_i^*)_{i\in F}$ is $C$-equivalent to the biorthogonals of $(e_i)_{i\in F}$. This is also true reversing the roles of $(e_i)$ and $(e_i^*)$. It is also not hard to see that if a $C$-Schreier unconditional sequence is not weakly null then it has a subsequences that generates an $\ell_1$ spreading model. Moreover, every infinite-dimensional Banach space  contains either $\ell_1$ or,  for every $\varepsilon>0$,  a $(2+\varepsilon)$-Schreier unconditional weakly null sequence  $ (e_i)$  \cite{O}.

The next result is a significant strengthening of \cite[Theorem 1.10]{DOS}.

\begin{prop}\label{prop: schreier}
Let $X$ be a Banach space with a $C$-Schreier unconditional Schauder basis $(e_i)_i$ so that all spreading models generated by subsequences of it are equivalent to a common sequence $(z_i)$, that is not equivalent to the unit vector basis of some $\ell_p$ or $c_0$. Then $|SP_{w}(X)|\vee|SP_{w}(X^*)| = \mathfrak{c}$.
\end{prop}

\begin{proof}
Since $(z_i)$ is not equivalent to the unit vector basis of $\ell_1$ then it must be 1-suppression unconditional and $(e_i)$ must be weakly null. In particular, it is subsymmetric. By the preceding discussion all spreading models generated by subsequences of $(e_i^*)$ are equivalent to $(z^*_i)_i$. By Theorem \ref{thm: l_p char2}, either $\overline{\spn(z_i)}$ or $\overline{\spn(z_i^*)}$ has a continuum of non-equivalent subsymmetric sequences generated by a vector. Let's assume the first is true as the second one is treated identically. By Proposition \ref{sm unc} there exists a continuum of block sequences $(u_i^\alpha)_{\alpha<\mathfrak{c}}$ that generate pairwise non-equivalent spreading models $(z^\alpha_i)$, all of which are unconditional. Therefore, at most one of them, say $(z^{\alpha_0}_i)$, is equivalent to the unit vector basis of $\ell_1$. This implies that $(u_i^\alpha)_{\alpha\neq\alpha_0}$ are all weakly null and therefore $|SP_{w}(X)| = \mathfrak{c}$.
\end{proof}

\begin{remark} If $X$ is as above,  then, combining Proposition~\ref{prop: schreier} and  Remark~\ref{rem: dominationorder},   either $SP_w(X)$ or $SP_w(X^*)$  contains a chain order-isomorphic to $(\mathbb{R}, \le)$  in  the domination order.  On the other hand,  it is known that if $Y$ is a separable Banach space and $SP_w(Y)$ is uncountable then $SP_w(Y)$ contains an antichain of cardinality $\mathfrak{c}$ in the domination order \cite{D, FR}.
\end{remark}

\end{document}